\documentclass{amsart}

\usepackage{amsfonts,amsmath,amsthm,latexsym,amscd,amsbsy,amssymb,fleqn,leqno,
pb-diagram,lamsarrow,graphicx}

\newtheorem{thm}{Theorem}[section]
\newtheorem{cor}[thm]{Corollary}
\newtheorem{lem}[thm]{Lemma}
\newtheorem{prop}[thm]{Proposition}
\theoremstyle{definition}

\theoremstyle{remark}

\numberwithin{equation}{section}

\newcommand{\invlim}{\displaystyle \lim _{\longleftarrow}}

\begin{document}

\title[Root closed function algebras on compacta of large dimension]
{Root closed function algebras on compacta of large dimension}

\author{N.~Brodskiy}
\address{University of Tennessee, Knoxville, TN 37996, USA}
\email{brodskiy@math.utk.edu}

\author{J.~Dydak}
\address{University of Tennessee, Knoxville, TN 37996, USA}
\email{dydak@math.utk.edu}

\author{A. Karasev}
\address{Nipissing University, North Bay,
ON, P1B 8L7, Canada}
\thanks{The third author was partially supported by NSERC Grant.}
\email{alexandk@nipissingu.ca}

\author{K. Kawamura}
\address{Institute of Mathematics, University of Tsukuba, Tsukuba, Ibaraki 305-8071, Japan}
\email{kawamura@math.tsukuba.as.jp}




\begin{abstract}
Let $X$ be a Hausdorff compact space and $C(X)$ be the algebra of
all continuous complex-valued functions on $X$, endowed with the
supremum norm. We say that $C(X)$ is (approximately) $n$-th root
closed if any function from $C(X)$ is (approximately) equal to the
$n$-th power of another function. We characterize the approximate
$n$-th root closedness of $C(X)$ in terms of $n$-divisibility of
first $\check {\rm C}$ech  cohomology groups of closed subsets of
$X$. Next, for each positive integer $m$ we construct
$m$-dimensional metrizable compactum $X$ such that $C(X)$ is
approximately $n$-th root closed for any $n$. Also, for each
positive integer $m$ we construct $m$-dimensional compact
Hausdorff space $X$ such that $C(X)$ is $n$-th root closed for any
$n$.
\end{abstract}

\maketitle


\section{Introduction}

Relations between algebraic closedness of the algebra of
continuous bounded complex-valued functions $C(X)$ on a space $X$
and topological properties of $X$ have been studied since 1960s
\cite{countryman}. Recall that the algebra $C(X)$ is called
algebraically closed if each monic polynomial with coefficients in
$C(X)$ has a root in $C(X)$. For locally connected compact
Hausdorff space, the algebra $C(X)$ is algebraically closed if and
only if ${\rm dim} X\le 1$ and $H ^1 (X;\mathbb{Z}) =0$
\cite{hatori_miura,miura_niijima}, where $H ^1 (X;\mathbb{Z})$
denotes the first {\v C}ech cohomology group of $X$ with the
integer coefficient (see section 2). It is proved in
\cite{miura_niijima} that for a first-countable compact Hausdorff
space $X$, algebraic closedness  of $C(X)$ is equivalent to a
weaker property of square root closedness. The latter means that
every function from $C(X)$ is a square of another function. It
should be noted that this property appears in the study of
subalgebras of $C(X)$~\cite{cirka}.

Even weaker property of approximate square root closedness was
introduced by Miura \cite{miura} and was proved to be equivalent to the
square root closedness when the underlying compact Hausdorff space $X$ is
locally connected.

There is a nice characterization of algebraic closedness of $C(X)$
when $X$ is a metrizable continuum. Namely, in this
case $C(X)$ is algebraically closed if and only if  $X$ is a
dendrite (i.e. a Peano continuum containing no simple closed
curves) \cite{kawamura_miura,miura_niijima}.

The approximate $n$-th root closedness of $C(X)$ was studied by
Kawamura and Miura and was proved to be equivalent to
$n$-divisibility of $H ^1(X;\mathbb{Z})$ under the additional
assumption ${\rm dim} X\le 1$. The universal space for metrizable
compacta with the approximately $n$-th root closed $C(X)$ is
constructed in \cite{ckkv}.

In this paper we characterize the approximate $n$-th root
closedness of $C(X)$ for any Hausdorff paracompact space $X$.
Namely, $C(X)$ is approximately $n$-th root closed if and only if
the group $H ^1 (A;\mathbb{Z})$ is $n$-divisible for every closed
subset $A$ of $X$. If ${\rm dim} X\le 1$, then the
$n$-divisibility of $H ^1 (X;\mathbb{Z})$ implies the
$n$-divisibility of $H ^1 (A;\mathbb{Z})$, so this generalizes
Theorem 1.3 of \cite{kawamura_miura}. Further, for each positive
integer $m$ we construct $m$-dimensional metrizable compactum $X$
such that $C(X)$ is approximately $n$-th root closed for any $n$.
Note that such examples were known in dimension 1 only. Also, for
each positive integer $m$ we construct $m$-dimensional compact
Hausdorff space $X$ such that $C(X)$ is $n$-th root closed for any
$n$. This example solves the problem posed in
\cite{kawamura_miura}: for a compact Hausdorff space $X$, does
square root closedness of $C(X)$ imply ${\rm dim} X\le 1$?

\section{Notations, definitions, and ideas of constructions}\label{preliminary}

All maps considered in this paper are continuous. For spaces $X$
and $Y$, we denote the set of all maps from $X$ to $Y$ by
$C(X,Y)$. As usual, by $\mathbb{Z}$, $\mathbb{Q}$, and
$\mathbb{C}$ we denote the integers, the rational numbers, and the
complex numbers, respectively. We let $\mathbb{C} ^*
=\mathbb{C}\setminus\{0\}$ be the multiplicative subgroup of
$\mathbb C$. An inverse spectrum over a directed partially ordered
set $(\mathcal A, <)$ consisting of spaces $X_{\alpha}$, $\alpha
\in\mathcal A$, and projections $p^{\beta} _{\alpha}\colon
X_{\beta}\to X_{\alpha}$, $\alpha ,\beta\in \mathcal A$, $\beta >
\alpha$, is denoted by $\{ X_{\alpha}, p^{\beta} _{\alpha},
\mathcal A\}$. Throughout this section $n> 1$ denotes an integer.

By $H ^k (X;G)$ we denote $k$-th $\check {\rm C}$ech cohomology
group of the space $X$ with abelian coefficient group $G$. Note
that in the case when $X$ is a Hausdorff paracompact space the
$\check {\rm C}$ech cohomologies are naturally isomorphic to the
Alexander-Spanier cohomologies \cite[p. 334]{spanier}. Note also
that due to the Huber's theorem \cite{huber} for a Hausdorff
paracompact space $X$ there exists a natural isomorphism between
the group of all homotopy classes of maps from $X$ to $K(G,k)$ and
the group $H ^k (X;G)$ if $G$ is countable. Here $K(G,k)$ denotes
the Eilenberg-MacLane complex.

For a space $X$, by $C(X)$ we denote the algebra of all bounded
complex-valued functions on $X$, endowed with the supremum norm.
We say that $C(X)$ is {\it approximately  $n$-th root closed} if for
every $f\in C(X)$ and every $\varepsilon > 0$ there exists $g\in
C(X)$ such that $||f-g^n|| < \varepsilon$. The algebra $C(X)$ is
said to be {\it $n$-th root closed } if any $f\in C(X)$ has an $n$-th root,
which means that there exists $g\in C(X)$ such that $f = g^n$.
Note that if $C(X)$ is (approximately) $n$-th root closed then
$C(A)$ is also (approximately) $n$-th root closed for any closed
subset $A$ of $X$.

We consider $C(X,\mathbb {C}^*)$ as a multiplicative subgroup of
$C(X)$ with metric inherited from $C(X)$. We say that $C(X,\mathbb
{C}^*)$ is (approximately) $n$-th root closed if any $f\in
C(X,\mathbb {C}^*)$ has an (approximate) $n$-th root in $C(X,\mathbb
{C}^*)$.

The basic idea explored in this paper --- the construction of
projective $n$-th root resolution --- is
outlined as follows.  The simplest case has been known in the theory of
uniform algebra and it is called the Cole construction
(cf.\cite[Chapter 3, \S 19, p.194-197]{stout}).

Given a space $X$ and a function $f\colon X\to\mathbb{C}$ it is
not always possible to solve the problem of finding an $n$-th root
of $f$ even approximately (consider for instance any homotopically
non-trivial map from a circle $S^1$ to $\mathbb{C}^*$).
Nevertheless, it is always possible to solve the $n$-th root problem
{\it projectively} in the following sense. There exists a space
denoted $R_n (X,f)$ and a map $\pi ^f\colon R_n (X,f) \to X$ such
that the composition $f\circ\pi ^f$ has an $n$-th root. The space
$R_n (X,f)$ is simply the graph of the (multivalued) $n$-th root
of $f$,
$$R_n(X,f) = \{ (x,z) \mid f(x) = z^n\}\subset X\times\mathbb{C},$$
and the map $\pi ^f$ is the natural projection on $X$. Obviously,
the projection of $R_n (X,f)$ to $\mathbb{C}$ is an $n$-th root of
the composition $f\circ\pi ^f$. We say that the space $R_n (X,f)$
together with the map $\pi ^f$ resolve the $n$-th root problem for
$f$ projectively.

Given any family of maps $\mathcal M\subset C(X)$ we can
projectively resolve the $n$-th root problem for all maps from
$\mathcal M$ simultaneously using the space
$$R_n (X,\mathcal M) = \{(x, (z_f)_{f\in \mathcal M})\mid f(x) = z^n_f\,\,\,\,\forall f\in
\mathcal M\}\subset X\times \mathbb{C} ^{\mathcal M}$$ and
defining $\pi ^{\mathcal M}\colon R_n (X,\mathcal M)\to X$ to be
the natural projection. Let $A$ and $B$ be two subsets of $C(X)$
such that $A\subset B$. There is a natural projection $\pi ^B _A
\colon R_{n}(X,B) \to R_{n}(X,A)$ defined by $\pi ^B _A [(x, (z_f)_{f\in
B})]= (x, (z_f)_{f\in A} )$. We let $R_n(X,\emptyset ) = X$ and
$\pi ^B _{\emptyset} = \pi ^{B}$.


Now we outline the ideas of our constructions in
Sections~\ref{approx roots} and~\ref{exact roots}. Suppose that we
want to construct a space $X$ with $n$-th root closed $C(X)$. Take
any space $X_1$ and resolve the $n$-th root problems for $X_1$
projectively using the space $X_2=R_n(X_1,C(X_1))$. Then resolve
all $n$-th root problems for $X_2$ projectively using $X_3$, and
so on. This way we obtain an inverse spectrum of spaces
$X_\lambda$ and define $X$ to be the inverse limit of this
spectrum. To guarantee that the $n$-th root problems for $X$ can
be solved, we need this spectrum to be {\it factorizing} in the
following sense: for any map $f\colon X\to\mathbb{C}$ there exist
a space $X_\lambda$ in the spectrum and a map $f_\lambda\colon
X_\lambda\to\mathbb{C}$ such that $f=f_\lambda\circ p_\lambda$
where $p_\lambda\colon X\to X_\lambda$ is the limit projection.
Then the projective resolution of the $n$-th root problem for $f_\lambda$
gives us a solution of the $n$-th root problem for $f$. In order to
obtain factorizing spectrum we make its length uncountable.
Namely, we construct the spectrum over $\omega _1$, the first
uncountable ordinal.

The space described above is not metrizable for two reasons.
First, the length of the spectrum used is not countable. Second,
for a metric compactum $X_\lambda$ and a subset $\mathcal M\subset
C(X_\lambda)$ the space $R_n (X_\lambda , \mathcal M)$ is
metrizable if and only if the set $\mathcal M$ is countable. If we
want $C(X)$ to have just the {\it approximate} $n$-th root
property, it is enough to construct a countable spectrum and for
each space $X_{\lambda}$ to resolve projective $n$-th root problem
for a countable {\it dense} set of maps from $C(X _{\lambda})$.
Then the limit space $X$ is a metrizable compactum, if we
start with a metrizable compactum $X_1$.

To guarantee that for the limit space $X =\invlim \{ X_{\lambda},
p^{\mu} _{\lambda}, \Lambda\}$ we can (approximately) solve the $n$-th
root problem for any function from $C(X)$ and for any $n>1$, we
represent the index set $\Lambda$ as the union of disjoint cofinal
subsets $\{\Lambda _n\} _{n=2}^{\infty}$. Then we construct the
spectrum by transfinite induction so that the space $X_{\lambda
+1}$ and the projection $p^{\lambda +1} _{\lambda}$ resolve
projectively (almost) all $n$-th root problems on $X_{\lambda}$,
where $\lambda\in\Lambda _{n}$. Since every set $\{\Lambda _n\}$
is cofinal, for any $n$ and any $\alpha$, (almost) every
$n$-th root problem on $X_\alpha$ will be projectively resolved at
some level $\lambda>\alpha$ where $\lambda\in\Lambda _{n}$.

To guarantee that the limit space $X$ has dimension ${\rm dim}
X\ge m$ we start the construction with the space $X_1$
homeomorphic to $m$-dimensional sphere $S^m$. Then we show that
the homomorphism $(p_1)^* \colon H^m (S^m;\mathbb{Q})\to H^m
(X;\mathbb{Q})$ induced by the limit projection is a monomorphism.
Therefore the mapping $p_1\colon X\to S^m$ is essential and hence
${\rm dim } X \ge m$.
To prove that the homomorphism above is a monomorphism we use a
construction called {\it transfer}, that briefly can be described
as follows. Suppose $G$ is a finite group acting on a compact
Hausdorff space $Y$. Let $Y/G$ be the quotient space and
$\pi\colon Y\to Y/G$ be the natural projection. Then there exists
a homomorphism $\mu ^* \colon H^* (Y;\mathbb{Q})\to H^*
(Y/G;\mathbb{Q})$ such that the composition $\mu ^*\pi ^*$ is the
multiplication by the order of $G$ in the group $H^*
(Y/G;\mathbb{Q})$. Therefore $\pi ^*$ is a monomorphism. See
Chapter II, $\S$19 of \cite{bredon} for more information on
transfer.

\section{Projective resolutions}

In this Section we establish some properties of projective
resolutions needed for our constructions in Sections~\ref{approx
roots} and~\ref{exact roots}. We begin with a summary of basic
properties of  space $R_n (X,\mathcal M)$.

\begin{prop}
 Let $X$ be a space, $\mathcal M$
be a subset of $C(X)$, and $n>1$ be an integer.

{\rm (a)} $R_n (X,\mathcal M)$ is the pull-back in the following
diagram:
$$
\begin{diagram}
\node{R_{n}(X,\mathcal M)}\arrow{e}\arrow{s,r} {\displaystyle \pi
^{\mathcal M}}\node{\mathbb{C} ^{\mathcal M}}
\arrow{s,r}{N}\\
\node{X}\arrow{e,t}{F}\node{\mathbb{C} ^{\mathcal M}}
\end{diagram}
$$
where $F\colon X \to \mathbb{C}$ is defined by $F(x)=
(f(x))_{f\in\mathcal M}$ and $N\colon \mathbb{C} ^{\mathcal  M}\to
\mathbb{C} ^{\mathcal M}$ is defined by $N((z_f)_{f\in \mathcal
M})= (z^n_f)_{f\in \mathcal M}$.

{\rm (b)} For any $f\in\mathcal M$ there exists $g\in
C(R_n(X,\mathcal M))$  such that $f\circ\pi ^{\mathcal M} = g^n$

{\rm (c)} If $X$ is a compact Hausdorff space then $R_n
(X,\mathcal M)$ is also a compact Hausdorff space and ${\rm dim}
R_n (X,\mathcal M)\le {\rm dim} X$.

\end{prop}

\begin{proof}
The statement (a) is obvious. To prove (b) we just let
$g[(x,(z_h)_{h\in A})] =z_f$.  To verify (c) we note first of all
that $R_n(X,\mathcal M)$ is a subset of the product   $X \times
\Pi _{f\in \mathcal M}\{z \mid z^n \in f(X)\}$ of compact
Hausdorff spaces. Moreover, $R_n (X,\mathcal M)$ is closed in this
product due to (a) and the compactness follows. For the dimension
part, observe that $\pi ^{\mathcal M}$ has zero-dimensional fibers
and apply \cite[Theorem 3.3.10]{engelking2}.
\end{proof}

In what follows we shall omit the index $n$ when this does not
cause ambiguities.


\begin{prop}\label{natural} For any space $X$ and any two subsets $A$ and $B$ of
$C(X)$ there exists a natural homeomorphism $h\colon
R(R(X,A),B\circ \pi ^A) \to R(X, A\cup B)$, where $B\circ\pi ^A =
\{ f\circ\pi ^A \mid f\in B\}$. This homeomorphism makes the
following diagram commutative
$$
\begin{diagram}
\node{R(R(X,A),B\circ \pi^{A})}\arrow{e,t}{\displaystyle h}\arrow{s,r}{\displaystyle \pi ^{B\circ \pi ^A}}\node{R(X, A\cup B)}\arrow{s,r}{\displaystyle \pi ^{A\cup B}}\\
\node{R(X,A)}\arrow{e,t}{\displaystyle \pi ^A}\node{X}
\end{diagram}
$$
\end{prop}

\begin{proof}
Note that both $R(X, A\cup B)$ and $R(R(X,A),B\circ \pi ^A)$ can
be viewed as subsets of $X\times\mathbb{C} ^A\times\mathbb{C}^B$. Namely,
$$R(X, A\cup B) = \{ (x, (z_f) _{f\in A}, (z_{g})_{g\in
B})\mid z ^n_{f} = f(x), z^n _{g} = g(x)\}$$
$$R(R(X,A),B\circ \pi ^A) = \{ (x, (z_f) _{f\in A},
(z_{g\circ\pi ^A})_{g\in B})\mid z^n _f = f(x), z ^n_{g\circ\pi
^A} = (g\circ\pi^A)[(x,(z_f)_{f\in A})]\}.$$ It remains to note
that these subsets coincide since $$(g\circ\pi ^A)[(x,(z_f)_{f\in
A})] = g (\pi ^A [(x,(z_f)_{f\in A})]) = g(x)$$ by the definition
of $\pi ^A$.
\end{proof}

\begin{prop}\label{spectrum}
Let $X$ be a compact Hausdorff space and $S$ be a subset of
$C(X)$. Let $\mathcal A$ be a family of subsets of $S$, partially
ordered by inclusion. Assume that $\mathcal A$ is a directed set
with respect to this order and that $\cup \mathcal A = S$. Then
$R(X,S)$ is naturally homeomorphic to the limit of the inverse
spectrum $\{R(X,A), \pi ^{B} _{A}, \mathcal A\}$.
\end{prop}

\begin{proof} Put $\Re =\invlim \{R(X,A), \pi ^{B} _{A}, \mathcal
A\}$. Define $h_A \colon R(X,S)\to R(X,A)$ for each $A\in\mathcal
A$ letting $h _A = \pi ^S _A$. The family of maps $\{ h_A \mid
A\in\mathcal A\}$ induces the limit map $h\colon R(X,S)\to\Re$
\cite[Proposition 1.2.13]{ch}. We claim that $h$ is a
homeomorphism. Since both $R(X,S)$ and $\Re$ are Hausdorff
compacta it is enough to check that $h$ is bijective. Since all
maps $\pi^S _A$ are surjective, $h$ is surjective by Theorem
3.2.14 in \cite{engelking1}. To verify the injectivity, it is enough,
for any two distinct points from $R(X,S)$, to find $A\in \mathcal
A$ such that the images of these two points under $h_A$ are
distinct. Let $y = (x,(z_f)_{f\in S})$ and $y' = (x',(z'_f)_{f\in
S})$ be two distinct points from $R(X,S)$. If $x\ne x'$, then any
$A\in\mathcal A$ will do. Otherwise there exists $f\in S$ such
that $z_f\ne z'_f$. Since $\cup\mathcal A = S$ there exists
$A\in\mathcal A$ such that $f\in A$ and one can easily see that $h
_A (y) \ne h_A (y')$.
\end{proof}

Later we use the following special case of Corollary~14.6 from
\cite{bredon}.

\begin{prop}\label{continuity}
Let $\mathcal S = \{ X_{\alpha}, p^{\beta} _{\alpha}, \mathcal
A\}$ be an inverse spectrum consisting of Hausdorff compact
spaces. Then there exists natural isomorphism
$$\lim _{\longrightarrow} H ^* (X_{\alpha};\mathbb{Q}) \cong H^*
(\invlim {\mathcal S}; \mathbb{Q}).$$
\end{prop}


\begin{prop}\label{monomorphism} Let $X$ be a compact Hausdorff space and $S$
be any subset of $C(X)$. Then for any integer $n>1$
$$(\pi ^S) ^* \colon H ^* (X;\mathbb{Q})\to H ^* (R_n(X,S); \mathbb{Q})$$
is a monomorphism.
\end{prop}

\begin{proof} (i) First, we prove the proposition for any space $X$ and
a set $S$ consisting of a single function $f$. There is an action
of $\mathbb{Z} _n$ on $R_n(X,f)$ whose orbit space is $X$, with
$\pi ^{f}$ being the quotient map. Namely, represent $\mathbb{Z}
_n$ as the group of $n$-th roots of $1$ and put $g\cdot
(x,z_f)=(x,g\cdot z_f)$. The proposition now follows from
Theorem~19.1 in~\cite{bredon}. Repeating the argument and applying
Proposition~\ref{natural} finitely many times, we see that the
proposition holds for every finite set $S$.

(ii) Finally, let S be any subset of $C(X)$. Let $S _{\rm fin}$
denote the set of all finite subsets of $S$, partially ordered by
inclusion. Proposition \ref{spectrum} implies that $R_n(X,S)$ is
the limit of the inverse spectrum $\{R_n(X,A),\pi ^B _A, S _{\rm
fin}\}$. We apply step (i) of this proof to conclude that
$\left(\pi ^B_A\right) ^* \colon H ^* (R_n(X,A); \mathbb{Q})\to H
^* (R_n(X,B); \mathbb{Q})$ is a monomorphism for all $A\subset B$
in $S _{\rm fin}$. An application of Proposition~\ref{continuity}
completes the proof.
\end{proof}

\section{Characterizations}\label{char}

\begin{lem}\label{homotopy}
If a map $f\colon X\to\mathbb{C}^*$ has an $n$-th root, then any map
$g\colon X\to \mathbb{C} ^*$ which is homotopic to $f$ also has an
$n$-th root.
\end{lem}

\begin{proof}
Apply the homotopy lifting property to the $n$-th degree covering
map $\mathbb{C}^*\to \mathbb{C}^*$, $z\mapsto z^n$.
\end{proof}

\begin{lem}\label{zero}
Let $X$ be a normal space. The following conditions are
equivalent.

\begin{itemize}

\item[(a)] $C(X)$ is approximately $n$-th root closed.

\item[(b)] $C(A,\mathbb{C} ^*)$ is approximately $n$-th root
closed for any closed subset $A$ of $X$.

\item[(c)] $C(A,\mathbb{C} ^*)$ is  $n$-th root closed for any
closed subset $A$ of $X$.
\end{itemize}
\end{lem}

\begin{proof} For a positive number $r$, let
$A(0,r) =\{ z\in\mathbb{C} \colon |z|\ge r\}$ and $B(0,r) =\{
z\in\mathbb{C} \colon |z|\le r\}$. Let $\rho _{\varepsilon} \colon
\mathbb{C} ^{*} \to A(0,\epsilon)$ be the radial retraction. Note
that $\rho_{\varepsilon}$ is homotopic to the identity map of
$\mathbb{C} ^*$.

${\rm (a)}\Rightarrow {\rm (b)}$ Take $\varepsilon >0$ and
consider a closed subset $A$ of $X$. Pick $f\in C(A,\mathbb{C}
^*)$ and put $h = \rho _{\varepsilon} \circ f$. Extend $h$ to a
function $F$ on $X$, applying the hypothesis (a) to find an $n$-th
root $g$, and restricting $g$ to $A$, we obtain a function
$g\colon A\to\mathbb{C}$ such that $||h-g^n|| <\varepsilon /2$.
This condition guarantees that $g\in C(A,\mathbb{C} ^*)$. It is
easy to see that $||f - g^n||<\varepsilon + \varepsilon /2 <
2\varepsilon$.

${\rm (b)}\Rightarrow {\rm (c)}$ Again, consider $f\in
C(A,\mathbb{C} ^*)$, where $A$ is a closed subset of $X$, and put
$h = \rho _{\varepsilon} \circ f$. Note that $h$ is homotopic to
$f$. Find $g\colon A\to\mathbb{C} ^*$ such that $||h-g^n||
<\varepsilon /2$. This condition guarantees that $g^n$ is
homotopic to $h$ and hence to $f$. An application of
Lemma~\ref{homotopy} completes the proof.

${\rm (c)}\Rightarrow {\rm (a)}$  Take $f\in C(X)$ and fix $\varepsilon >0$.
Consider $A=f^{-1}(A(0,\varepsilon ))$ and
$B=f^{-1}(B(0,\varepsilon ))$. Find $g\in C(A,\mathbb{C}^*)$ such
 that $f|_A = g^n$. Note that $g(A\cap B)
\subset B(0,\sqrt[n]{\varepsilon})$ and we can extend $g$ over $X$
to $\overline{g}$ such that $\overline{g} (B) \subset
B(0,\sqrt[n]{\varepsilon})$. It is easy to check that $||f-
\overline{g} ^{\,n}||<2\varepsilon$.
\end{proof}

We let $S^1 = \{ z\in\mathbb{C}\colon |z|=1\}$. Suppose $Y$ is a
Hausdorff paracompact space. Huber's Theorem \cite{huber} implies
the existence of a canonical isomorphism $H^1 (Y;\mathbb{Z})\cong
[Y,S^1]$. Here $[Y,S^1]$ denotes the group of all homotopy classes
of maps from $Y$ to $S^1$ with the group operation induced by the
multiplication of maps in $C(Y,S^1)$. We denote the homotopy class
of a map $f\in C(Y,S^1)$ by $[f]$.

\begin{thm} Let $X$ be a Hausdorff paracompact space. Then $C(X)$
is approximately $n$-th root closed iff $H^1(A;\mathbb{Z})$ is
$n$-divisible for every closed subset $A$ of $X$.
\end{thm}

\begin{proof}
Consider a closed subset $A$ of $X$. First, suppose that $C(X)$ is
approximately $n$-th root closed.  Let $f\colon A\to S^1$ be a
representative of an arbitrary element of $H^1 (A;\mathbb{Z})$. By
condition (c) of Lemma~\ref{zero} there exist $g\colon A\to S^1$
such that $g^n =f$ and hence $n[g] = [f]$ in $H^1(A;\mathbb{Z})$.

In order to prove the converse part, we verify the condition (c)
of Lemma~\ref{zero}. Pick $f\in C(A,\mathbb{C}^*)$. Then $f$ is
homotopic to a map $\widetilde{f}\colon A\to S^1$. Since
$[\widetilde{f}\,\,]\in H^1(A;\mathbb{Z})$ is divisible by $n$
there exists $h\colon A\to S^1$ such that $h^n$ is homotopic to
$\widetilde{f}$ and hence to $f$. Lemma~\ref{homotopy} implies
that $f$ has an $n$-th root.
\end{proof}

\section{Compacta with approximately root closed $C(X)$}\label{approx roots}

\begin{lem}\label{invlim}
Let $\mathcal S = \{X_i,p^{i+1}_i\}$ be an inverse sequence of
compact metrizable spaces and let $X = \invlim\mathcal S$. Consider the
following two conditions:

\begin{itemize}

\item[(a)] $C(X)$ is approximately $n$-th root closed.

\item[(b)] For any
$i$, any closed subset $A_i$ of $X_i$ and any map $h\colon
A_i\to\mathbb{C} ^*$  there exists $j>i$ such that the map $h\circ
p^j_i\colon A_j\to\mathbb{C}^*$ has $n$-th root, where
$A_j=(p^j_i)^{-1}(A_i)$.
\end{itemize}

The condition (b) implies the condition (a).  Moreover if all projections of
$\mathcal S$ are
surjective, then the converse implication (a)$\rightarrow$(b) also holds.

\end{lem}


\begin{proof} Put $X=\invlim\mathcal S$. First, we show that
$C(X)$ is approximately $n$-th root closed by checking the
condition (b) of Lemma~\ref{zero}. Let $A$ be a closed subset of
$X$ and $f\in C(A,\mathbb{C} ^*)$ be a function. Take any
$\varepsilon>0$. There exist $i$ and a mapping $f_i\colon p_i (A)
\to \mathbb{C}^*$ such that $f_i\circ p_i |_A$ is
$\varepsilon$-close to $f$. Let $A_{i} = p_{i}(A)$ and find $j>i$
such that the map $f_i\circ p^j_i\colon A_j\to\mathbb{C}^*$ has an
$n$-th root.
Then $g = h \circ f_{j}$ is an $n$-th root of $f_i\circ p_i | _A$.
Obviously, $||f-g^n||<\varepsilon$.

Conversely, suppose $C(X)$ is approximately $n$-th root closed and
all projections of $\mathcal S$ are surjective. Pick $i$ and
consider a closed subset $A_i$ of $X_i$ and a map $h\colon
A_i\to\mathbb{C} ^*$. Let $\varepsilon = \min\{|h(x)|\colon x\in
A_i\}$. Put $A = (p_i)^{-1} (A_i)$. There exists $g\colon A\to
\mathbb{C}^*$ such that $g^n$ is $\varepsilon /4$-close to $h\circ
p_i|_A$. We can find $j>i$ and a map $g_j\colon p_j(A)
\to\mathbb{C} ^*$ such that $(g_j\circ p_j)^n$ is $\varepsilon
/4$-close to $g^n$. Let $A_j=(p^j_i)^{-1}(A_i)$.  Since all
projections of $\mathcal S$ are surjective, $p_j (A) = A_j$. Using this, it
is not hard to verify that $(g_j)^n$ is $\varepsilon /2$-close, and
hence homotopic, to $h\circ p^j_i$. Lemma \ref{homotopy} implies
that $h\circ p^j_i$ has an $n$-th root.
\end{proof}

\begin{thm} For every positive integer $m$ there exists an $m$-dimensional
compact metrizable space $X$ such that $C(X)$ is approximately
$n$-th root closed for all positive integers $n$.
\end{thm}

\begin{proof} We obtain $X$ as the inverse limit of a sequence
$\mathcal S = \{X_i,p^{i+1}_i \}$, consisting of $m$-dimensional
metrizable compacta. The sequence is constructed by induction as
follows. Represent the set of all positive integers as a union of
disjoint infinite subsets $\{\Lambda _n\} _{n=2}^{\infty}$. Put
$X_1=S^m$, the $m$-dimensional sphere. Suppose the space $X_k$ has
been already constructed. Fix a countable collection $\mathcal
B_k$ of closed subsets of $X_k$ such that for each closed subset
$A$ of $X_k$ and for any open neighborhood $U$ of $A$ there exists
$B \in \mathcal B_k$ such that $A \subset B \subset U$.
For each $B\in \mathcal B_k$ fix a family $\mathcal F ^{*} _{B}$
of maps from $B$ to $\mathbb{C} ^*$ which is dense in the space $C
(B,\mathbb{C}^*)$. For every map from the family $\mathcal
F^{*}_{B}$ we fix its extension to a map from $X_k$ to $\mathbb{C}$
and denote the family of these extensions by $\mathcal F _B$. Let
$\Phi _k = \cup\{\mathcal F_B \mid B\in\mathcal B _k\}$. Define
$X_{k+1} = R_n (X_k, \Phi _k)$ where $n$ is such that $k\in\Lambda
_n$, and let $p^{k+1}_k = \pi ^{\Phi _k}$.

Put $X =\invlim \mathcal S$. To verify that $C(X)$ is
approximately $n$-th root closed for each $n >1$ it is enough to
show that condition (b) of Lemma \ref{invlim} is satisfied for the
inverse sequence $\mathcal S$. Fix $n>1$. Pick $i$ and consider a
closed subset $A_i$ of $X_i$ and a function $h\colon
A_i\to\mathbb{C}^*$. Take a number $j>i$ such that $j-1 \in\Lambda
_n$. Let $A_j =\left(p^j_i\right)^{-1} (A_i)$. We show that the
map $h\circ p^j_i\colon A_j \to\mathbb{C} ^*$ has an $n$-th root.
Put $A_{j-1} =\left(p^{j-1}_i\right)^{-1} (A_i)$. Let $g$ be an
extension of the map $h\circ p^{j-1}_i\colon A_{j-1} \to\mathbb{C}
^*$ to some neighborhood $U$ of $A_{j-1}$. There exists
$B\in\mathcal B _k$ and a function $f\colon B\to\mathbb{C}^*$ such
that $A_{j-1}\subset B\subset U$ and the restriction $g|_B$ is
homotopic to $f$. Let $\widetilde{f}\colon B\to\mathbb{C}$ be the
extension of $f$ that belongs to the family $\Phi _k$. Since the
map $p^j_{j-1}$ resolves the projective $n$-th root problem for
$\widetilde{f}$ the map $f\circ p^j_{j-1}|_{A_j}$ has an $n$-th
root. By Lemma~\ref{homotopy} the map $g\circ p^j_{j-1}|_{A_j}$
has an $n$-th root. It remains to note that $h\circ p^j_i |_{A_j}
= g\circ p^j_{j-1}|_{A_j}$.

 Note that ${\rm dim} X \le m$ since
all $X_k$ are at most $m$-dimensional.
Proposition~\ref{monomorphism} implies that $(p^{k+1}_k)^*\colon
H^m (X_k;\mathbb{Q})\to H^m (X_{k+1};\mathbb{Q})$ is a
monomorphism. Applying Proposition~\ref{continuity} we conclude
that $(p_1)^* \colon H^m (S^m;\mathbb{Q})\to H^m (X;\mathbb{Q})$
is a monomorphism. Thus the limit projection $p_1\colon X\to S^m$
is essential and therefore ${\rm dim} X \ge m$.
\end{proof}

Let $\mathcal K$ be a class of spaces.  A space $Z \in \mathcal K$
is called a universal space for the class $\mathcal K$, if every
space in $\mathcal K$ is topologically embedded in $Z$. For an
positive integer $n$ and $\tau \geq \omega$, let ${\mathcal
A}_{\tau}(n)$ (${\mathcal A}_{\tau}$ resp.) be the class of all
compact Hausdorff spaces $X$ such that $w(X) \leq \tau$ and $C(X)$
is $n$-th root closed ($C(X)$ is $n$-th root closed for each $n >
1$ resp.). It was shown in \cite{ckkv}, Corollary 1.3 that
${\mathcal A}_{\tau}(n)$ contains a univesal space for any $\tau
\geq \omega$ and any $n>1$.  Using the idea of the proof of
Theorem 1.2 from \cite{ckkv} one can show that ${\mathcal
A}_{\tau}$ also contains a universal space.

\begin{cor} Let $Y$ be a universal space with respect to the class
${\mathcal A}_{\omega}$ or ${\mathcal A}_{\omega}(n)$. Then $Y$ is
infinite-dimensional.
\end{cor}

Hence, any universal space for the
class ${\mathcal A}_{\tau}(n)$ (${\mathcal A}_{\tau}$ resp.) must be
infinite dimensional for any $\tau \geq \omega$.

Also we may consider the subclass ${\mathcal A}_{m, \tau}(n)$
(${\mathcal A}_{m, \tau}$ resp.) consisting of all spaces in
${\mathcal A}_{\tau}(n)$ (${\mathcal A}_{\tau}$ resp.) of
dimension at most $m$.  Theorem 1.2 of \cite{ckkv} also proves
that the class ${\mathcal A}_{1, \tau}(n)$ contains a universal
space.  A similar proof, based on the Marde{\v s}i{\v c}
factorization theorem \cite{m}, works to prove that the class
${\mathcal A}_{m, \tau}(n)$ (${\mathcal A}_{m, \tau}$ resp.)
contains a universal space.

\section{Compacta with root closed $C(X)$}\label{exact roots}

In this section, for any positive integer $m$ we construct
a compact Hausdorff space $X$ with  ${\rm dim} X = m$ such that
$C(X)$ is $n$-th root closed for all $n$. Note that for a
metrizable continuum $Y$ the algebra $C(Y)$ is square root closed
if and only if $Y$ is a dendrite, and therefore ${\rm dim}\, Y\le
1$ \cite{kawamura_miura,miura_niijima}. This forces the space $X$
above to be non-metrizable.

\begin{lem}\label{factorizing}
Let $\mathcal S =\{ X_{\alpha}, p^{\beta}_{\alpha},\mathcal A\}$
be a factorizing spectrum. In order for $C(\invlim\mathcal S)$ to
be $n$-th root closed it is sufficient that for any $\alpha\in
\mathcal A$ and any function $h\in C(X_{\alpha})$ there exists
$\beta >\alpha$ such that $h\circ p^{\beta}_{\alpha}$ has an
$n$-th root.
If all limit projections of $\mathcal S$ are surjective the above
condition is also necessary.
\end{lem}
\begin{proof} Put $X =\invlim\mathcal S$. Consider $f\in C(X)$.
Since $\mathcal S$ is factorizing there exists $\alpha$ and
$f_{\alpha}\in C(X_{\alpha})$ such that $f = f_{\alpha}\circ
p_{\alpha}$. By the condition of the lemma we can find $\beta
>\alpha$ and $g_{\beta} \colon X_{\beta}\to\mathbb{C}$ such that
$(g_{\beta})^n = f_{\alpha}\circ p^{\beta}_{\alpha}$. It is easy
to verify that $g = g_{\beta}\circ p_{\beta}$ is $n$-th root of
$f$.

 Now suppose that all limit
projections of $\mathcal S$ are surjective and $C(X)$ is $n$-th
root closed. Consider $\alpha\in \mathcal A$ and $h\in
C(X_{\alpha})$. There exists $g\in C(X)$ such that $g ^n = h\circ
p_{\alpha}$. Since $\mathcal S$ is factorizing there exists $\beta
>\alpha$ and $g_{\beta}\colon X_{\beta} \to \mathbb{C}$ such that
$g = g_{\beta}\circ p_{\beta}$. Since the projection $p_{\beta}$
is surjective $(g_{\beta})^n = h\circ p^{\beta}_{\alpha}$.
\end{proof}

\begin{thm}\label{example_square-root}
For each positive integer $m$, there exists a compact Hausdorff
space $X$  with ${\rm dim} X = m$  and such that $C(X)$ is $n$-th
root closed for any $n$.
\end{thm}

\begin{proof}
Represent the ordinal $\omega _1$ as the union of countably many
of disjoint uncountable subsets $\{ \Lambda _n \}_{n=2}^{\infty}$.
Starting with $X_0 = S^m$, where $S^m$ denotes $m$-dimensional
sphere, by transfinite induction we define an inverse spectrum
$\mathcal S = \{ X_{\alpha}, p ^{\beta}_{\alpha},\omega _1\}$ as
follows. If $\beta =\alpha +1$ then define $X_{\beta} = R_n
(X_{\alpha},C(X_{\alpha}))$, where $n$ is such that $\alpha\in
\Lambda _n$, and let $p ^{\beta} _{\alpha} =\pi ^{C(X_{\alpha})}$.
If $\beta$ is a limit ordinal, then define $\displaystyle
X_{\beta} =\lim_{\longleftarrow} \{ X_{\alpha},
p^{\gamma}_{\alpha}, \alpha <\beta\}$ and, for $\alpha < \beta$,
let $\displaystyle p ^{\beta}_{\alpha}$ be the limit projection.

Put $X =\invlim \mathcal S$. To verify that $C(X)$ is
 $n$-th root closed for each $n >1$ it is enough to check the condition
 of Lemma \ref{factorizing} for the spectrum
 $\mathcal S$. Consider $n>1$. Since the spectrum $\mathcal S$ has
 length $\omega _1$ it is factorizing \cite[Corllary 1.3.2]{ch}.
Consider a function $h\colon X_{\alpha}\to\mathbb{C}$ and take an
ordinal $\gamma >\alpha$ such that $\gamma\in \Lambda _n$. Since
the map $p^{\gamma +1} _{\gamma}$ resolves projective $n$-th root
problem for $h\circ p^{\gamma}_{\alpha}$, the map $h\circ
p^{\gamma +1}_{\alpha}$ has an $n$-th root.

Note that ${\rm dim} X_{\alpha} \le m$ for each $\alpha$ and hence
${\rm dim} X \le m$. We claim that $(p _{\alpha} ^{\beta} )^{*}
\colon H ^* (X _{\alpha};\mathbb{Q})\to H ^* (X
_{\alpha};\mathbb{Q})$ is a monomorphism for all $\alpha < \beta <
\omega _1$. Indeed, in the case $\beta = \alpha +1$ it follows
from Proposition \ref{monomorphism}, and then in general case it
is due to Proposition~\ref{continuity}. Finally, again with the
help of Proposition~\ref{continuity}, we conclude that $p _0 ^*
\colon H ^m (S^m;\mathbb{Q}) \to H ^m (X;\mathbb{Q})$ is a
monomorphism and hence the map $p_0 \colon X\to S^m$ is essential.
This implies ${\rm dim} X \ge m$.
\end{proof}

It is not hard to verify that if $C(Y)$ is $n$-th root closed for
some (completely regular) space $Y$ then $C(\beta Y)$ is also
$n$-th root closed. Here by $\beta Y$ we denote the Stone-$\check
{\rm C}$ech compactification of $Y$.

\begin{cor}
There exists an infinite dimensional compact Hausdorff space $X$
such that $C(X)$ is $n$-th root closed for all $n$.
\end{cor}
\begin{proof}
For each $m$, let $X_m$ denote compactum provided by Theorem
\ref{example_square-root}. We put $X = \beta (\oplus \{X_m \mid
m\in \omega \})$.
\end{proof}




\end{document}